
\input amssym.def
\input amssym.tex
\input xy
\input verbatim
\xyoption{all}

\def\cDer{{\cal D}\!er}

\def\cHom{{\cal H}\!o\!m}
\def\cOmega{{\it\Omega}}
\def\cTheta{{\it\Theta}}
\def\cDiff{{\cal D}\!\hbox{\it iff}}

\def\cO{{\cal O}}
\def\cD{{\cal D}}
\def\cE{{\cal E}}

\def\cH{{\cal H}}

\def\cL{{\cal L}}
\def\cM{{\cal M}}
\def\cN{{\cal N}}

\def\cP{{\cal P}}

\def\bL{{\bf L}}
\def\bR{{\bf R}}
\def\Gr{{\rm Gr}}

\def\Im{{\rm Im}}

\def\lMod{{\rm\hbox{-}Mod}}
\def\rMod{{\rm Mod\hbox{-}}}
\def\bMod{{\rm\hbox{-}Mod\hbox{-}}}
\def\lDiff{{\rm\hbox{-}Diff}}

\def\DR{{\rm DR}}
\def\tDR{{\widetilde{\rm DR}}}
\def\Ler{{\rm Ler}}
\def\point{{\scriptscriptstyle\bullet}}
\def\ot{\leftarrow}
\def\longmapsot{\hskip-2pt\longleftarrow\hskip-1.2pt\mapstochar}
\def\isom{\cong}
\def\id{{\rm id}}
\def\tot{{\rm tot}}

\centerline{\bf Algebraic Connections vs. Algebraic {$\cD$}-modules:
inverse and direct images.}
\medskip
\centerline{ by Maurizio Cailotto and Luisa Fiorot}
\medskip
\centerline{ Department of Mathematics - University of Padova (Italy)}
\medskip
\centerline{\tt \verbatim9  [cailotto,fiorot]@math.unipd.it 9}

\bigskip
{\bf Introduction}
\medskip
In the dictionary between the language of
(algebraic integrable) connections and that of (algebraic) $\cD$-modules,
the operations of direct and inverse images for a smooth morphism
are very important.
To compare the definitions of inverse images for connections and
$\cD$-modules
is easy.
But the comparison between direct images for connections
(the classical construction of the Gauss-Manin connection for smooth
morphisms)
and for $\cD$-modules, although known to specialists,
has been explicitly proved only recently in a paper of
Dimca, Maaref, Sabbah and Saito in 2000 (see [DMMS]),
where the authors' main technical tool was M. Saito's equivalence between
the derived category of $\cD$-modules and a localized category of
differential complexes.
\endgraf
The aim of this short paper is to give a simplified  summary of the [DMSS]
argument,
and to propose an alternative proof of this comparison
which is simpler, in the sense that it does not use  Saito equivalence.
Moreover, our alternative strategy of comparison works in a context which
is a precursor
to the Gauss-Manin connection
(at the level of $f^{-1}\cD_Y$-modules, for a morphism $f:X\to Y$),
and may be of some intrinsic interest.
 \endgraf
In section 1 we recall some generalities on connections and $\cD$-modules.
In section 2 we compare the operations of ``inverse image''
for connections and $\cD$-modules.
Section 3 is devoted to the comparison of
the  Gauss-Manin connection (in the case of smooth morphisms)
with the notion  of direct image for $\cD$-modules: we supply a
simplified summary of the [DMSS] argument.
Finally, in the last section we propose our alternative proof of this
comparison
which does not use  Saito equivalence.

\bigskip
{\bf\S 1. Generalities on connections and {$\cD$}-modules}
\medskip
Let $X$ be a smooth $K$-variety, where $K$ is a field of characteristic
$0$.
Following the terminology of [EGA$I\!V$, \S 16]
we denote by $\cOmega^1_X$ the ${\cO}_X$-module of differentials,
by $\cP_X^1$ the ${\cO}_X$-algebra of principal parts of order one:
its two structures as ${\cO}_X$-algebra will be referred to as the 
``left'' and ``right'' structures, and  tensor products will be specified by the
position of the $\cP_X^1$ factor.
We also use $\cDer_X$ or $\cTheta_X$ to denote the ${\cO}_X$-module of derivations
(${\cO}_X$-dual of $\cOmega^1_X$, endowed with the usual structure of
Lie-algebra),
and  ${\cD}_X$ to indicate the graded (left) ${\cO}_X$-algebra of differential
operators.
On ${\cD}_X$ we consider the increasing filtration $F$ defined by
the order of differential operators.
Then the associated graded ${\cO}_X$-algebra, denoted by $\Gr\cD_X$,
is commutative and it is generated (as ${\cO}_X$-algebra) by
$\cDer_X\subseteq F^1{\cD}_X$.
\endgraf
For any ${\cO}_X$-module ${\cE}$ we will use the standard notation
$\cP_X^1({\cE})$ for $\cP_X^1\otimes_{{\cO}_X}{\cE}$,
where the tensor product involves the right ${\cO}_X$-module structure
of $\cP_X^1$, while  the ${\cO}_X$-module structure is
given by the left  ${\cO}_X$-module structure  on $\cP_X^1$.

\medskip
{\bf 1.1 Connections and $\cD$-modules}
\medskip
\noindent 
Let $\cE$ be an ${\cO}_X$-module.
The following supplementary structures on $\cE$ are equivalent:
\itemitem{\bf (i)}
an integrable connection, that is a morphism of abelian sheaves
$\nabla:{\cE}\to\cOmega^1_X\otimes_{{\cO}_X}{\cE}$
Leibniz w.r.t. sections of ${\cO}_X$ and such that $\nabla^2=0$
for the natural extension of $\nabla$ to the De Rham sequence;
\itemitem{\bf (ii)}
an $\cO_X$-linear section $\delta:{\cE}\to\cP^1_X\otimes_{{\cO}_X}{\cE}$
of the canonical morphism $\pi:\cP^1_X\otimes_{{\cO}_X}{\cE}\to\cE$;
\itemitem{\bf (iii)}
an  $\cO_X$-linear Lie-algebra homomorphism
$\Delta:\cDer_X\to\cDiff_{X}({\cE})$ (for the usual Lie-algebra
structures), where $\cDiff_{X}({\cE})$ is the sheaf of differential
operators of $\cE$;
\itemitem{\bf (iv)}
a structure of left algebraic $\cD_X$-module on $\cE$.
\endgraf
The dictionary between these equivalent structures is
well explained in [BO, 2.9, 2.11, 2.15]; let us give a sketch.
\par
If $c=c_X({\cE}):{\cE}\to\cP^1_X\otimes_{{\cO}_X}{\cE}$ denotes the
canonical
inclusion,
then $\delta=c+\nabla$ and $\nabla=\delta-c$.
\par
For  any $\partial$ section of $\cDer_X$
the morphism $\Delta$ is defined by
$\Delta_{\partial}=(\partial\otimes\id)\circ\nabla$,
i.e. $\Delta_{\partial}(e)=\langle\partial,\nabla(e)\rangle$.
On the other hand,  the reconstruction of $\nabla$ from $\Delta$ requires a
description
using local coordinates $x_i$ on $X$ 
($dx_i$ and $\partial_i$ are the dual bases of differentials and derivations):
if $e$ is a section of $\cE$, then
$\nabla(e)=\sum_idx_i\otimes\Delta_{\partial_i}(e)$.
\par
The morphism $\Delta$ is equivalent to the data of a
left $\cD_X$-module structure on  $\cE$ since it extends
to a left action of $\cD_X$ on $\cE$ ([Bo,VI,1.6]).

\endgraf
\medskip
\noindent{\bf 1.2}
A morphism of connections on $X$ is an ${\cO}_X$-linear morphism
$h:{\cE}\to{\cE}'$ compatible with the data, that is, such that
$\nabla'\circ h=({\rm id}\otimes h)\circ\nabla$,
or $\delta'\circ h=({\rm id}\otimes h)\circ\delta$,
or equivalently
$\Delta'_\partial\circ h=h\circ\Delta_\partial$
for any $\partial$ section of $\cDer_X$,
or finally which is $\cD_X$-linear.
\endgraf

\bigskip
{\bf\S 2. Inverse image for connections and {$\cD$}-modules}
\medskip
Let $f:X\to Y$ be a finite type morphism of smooth $K$-varieties.
For any ${\cO}_Y$-module ${\cE}$, let
$f^\ast(\cE)=\cO_X\otimes_{f^{-1}\cO_Y} f^{-1}\cE$
the inverse image by $f$.
\medskip
{\bf 2.1 Inverse image for connections}
\medskip
The easiest definition for the inverse image by $f$ of a connection $\cE$
on $Y$
is given in terms of $\cO_Y$-linear maps.
If $\delta:\cE\to\cP^1_Y\otimes_{\cO_Y}\cE$ defines the connection,
 let  $f^\ast\delta$ be the composition of the inverse image of
$\delta$ with
the canonical morphism $f^\ast(\cP^1_Y)\to\cP^1_X$.
Then we have a morphism
$$
f^\ast\delta: f^\ast\cE\to\cP^1_X\otimes_{\cO_X} f^\ast\cE
\isom
\cP^1_X\otimes_{f^{-1}\cO_X} f^{-1}\cE\ ,
$$
which is clearly an $\cO_X$-linear section of the canonical map
$\pi:\cP^1_X\otimes_{\cO_X}f^\ast\cE\to f^\ast\cE$.
\endgraf
An explicit description of the connection $f^\ast\nabla$
on $f^\ast\cE$ can be given in the following way:
$$
\eqalign{
(f^\ast\nabla)(\alpha\otimes e)
&=(f^\ast\delta-c_X)(\alpha\otimes e)\cr
&= \alpha f^\ast\delta(1\otimes e)-c_X(\alpha\otimes e) \cr
&= \alpha f^\ast((\nabla+c_Y)(e))-{\Bbb I}\otimes(\alpha\otimes e) \cr
&= \alpha f^\ast(\nabla(e))+\alpha({\Bbb I}\otimes e)-{\Bbb
I}\otimes(\alpha\otimes e) \cr
&= \alpha \nabla(e)+\alpha\otimes 1\otimes e-1\otimes\alpha\otimes e \cr
&= \alpha \nabla(e)+d(\alpha)\otimes e \cr
}
$$
(as usual, $\alpha$ is a section of $\cO_X$
and $e$ is a section of $\cE$, or $f^{-1}\cE$).
\medskip
{\bf 2.2 Inverse image for $\cD$-modules}
\medskip
Let $\cM$ be a left
$\cD_Y$-module. The inverse image as $\cO$-modules
$f^\ast\cM=\cO_X\otimes_{f^{-1}\cO_Y}f^{-1}\cM$ locally admits an action
of $\cD_X$ defined by
$$
\alpha'(\alpha\otimes m)=(\alpha'\alpha)\otimes m
\quad
\hbox{and}
\quad
\partial(\alpha\otimes m)=\partial(\alpha)\otimes m +
\alpha\big(\sum\nolimits_i\partial(y_i)\otimes\eta_i(m) \big)
$$
where $\partial$ is a section of $\cD_X$,
$m$ a section of $\cM$ (or $f^{-1}\cM$), $\alpha,\alpha'$ sections of
$\cO_X$
($y_i$ local coordinates on $Y$ and $\eta_i$ the dual basis of
$dy_i$).
These local definitions globalize to a
$\cD_X$-module structure
on $f^\ast\cM$ [Bo,VI.4].
\endgraf
In this way the $\cO_X$-module
$f^\ast\cD_Y=\cO_X\otimes_{f^{-1}\cO_Y}f^{-1}\cD_Y$
is endowed  with a structure of left $\cD_X$-module,
compatible with the obvious structure of  $f^{-1}\cD_Y$-module
(by right multiplication). With this structure, $f^\ast\cD_Y$ is usually
denoted by $\cD_{X{\to}Y}$
and the inverse image of a left $\cD_Y$-module $\cM$
can be defined as
$$
f^\ast\cM=\cD_{X{\to}Y}\otimes_{f^{-1}\cD_Y}f^{-1}\cM\ ,
$$
taking account also of the structure as $\cD_X$-module.
\medskip
{\bf 2.3 Comparison}
\medskip
Let $\cM$ be a $\cD_Y$-module. We regard it as a connection on $Y$
and  consider its inverse image as a connection.
The action of derivations
is described in terms of local coordinates $y_i$ on $Y$, by
$$
\eqalign{
(f^\ast\Delta)_\partial(\alpha\otimes m)
&= \langle\partial, (f^\ast\nabla)(\alpha\otimes m) \rangle\cr
&= \langle\partial, \alpha \nabla(m)+d(\alpha)\otimes m \rangle\cr
&= \langle\partial,
\alpha(\sum\nolimits_i dy_i\otimes\Delta_{\eta_i}(m))+d(\alpha)\otimes
m\rangle\cr
&=\alpha\big(\sum\nolimits_i\partial(y_i)\otimes\Delta_{\eta_i}(m) \big)
    +\partial(\alpha)\otimes m
}
$$
where $\partial$ is a section of $\cDer_X$,
$m$ is a section of $\cM$ (or $f^{-1}\cM$), $\alpha$ is a section of
$\cO_X$
(and $\eta_i$ is the dual basis of $dy_i$).
Thefore the local descriptions
make clear that for a connection $\cE$ on $Y$,
the inverse image as a connection induces the structure of $\cD_X$-module
given by the inverse image of the corresponding $\cD_Y$-module.

\bigskip
{\bf\S 3 Direct image for connections and {$\cD$}-modules (and comparison
following [DMSS])}
\bigskip
Let $f:X\to Y$ be a smooth morphism of smooth $K$-varieties.
In order to compare the notions of (derived) direct images in the category
of connections (the Gauss-Manin connections) and in the category of
$\cD$-modules,
we need some preliminary materials, most concerning right
$\cD$-modules, De Rham functors,
differential complexes (and the M.Saito equivalence). 
In this paper, we use the convention that the shift of complexes 
does not change the sign of differentials. 
\smallskip
{\bf 3.1 Right and left $\cD$-modules}
\smallskip
We denote by $\cD\lMod$ the category of left $\cD$-modules
and by $\rMod\cD$ the category of right $\cD$-modules.
It is well known that $\omega_X=\cOmega_X^{\dim X}$
has a canonical structure of right $\cD_X$-module (see [Bo,VI,3.2]).
Let us define $\omega_X(\cD_X)=\omega_X\otimes_{\cO_X}\cD_X$.
It is endowed with two different structures of right $\cD_X$-module:
the first comes from the right multiplication on $\cD_X$ and the other
is induced by the tensor product over $\cO_X$ of a right and a left
$\cD_X$-module.
There exists a unique involution 
$\iota:\omega_X(\cD_X)\rightarrow\omega_X(\cD_X)$
which is the identity on $\omega_X$ and exchanges these two right
$\cD_X$-module structures (see [Sa,1.7], 
using local coordinates $x_i$ on $X$, 
the involution sends $\omega\otimes P$ to $\omega\otimes P^\ast$, 
where $\omega=dx_1\wedge\cdots\wedge dx_{d_X}$, and 
$P^\ast$ is the transposition of $P$, defined by 
$\alpha^\ast=\alpha$ for sections of $\cO_X$, 
$\partial_i^\ast=-\partial_i$, and $(PQ)^\ast=Q^\ast P^\ast$).
In the same way, we define
$\omega_X^{-1}(\cD_X)=\cD_X\otimes_{\cO_X}\omega_X^{-1}
=\cHom_{\cD_X}(\omega_X(\cD_X),\cD_X)$ and
we notice that $\omega_X^{-1}(\cD_X)$ has two compatible and
``interchangeable''
structures of left $\cD_X$-module.
\endgraf
We have an equivalence of categories between $\cD_X\lMod$ and $\rMod\cD_X$
given by the quasi inverse functors:
$$
\matrix{
\cD_X\lMod &\longleftrightarrow &\rMod\cD_X \cr
\hfill\cM    &\longmapsto
&\omega_X(\cM)=\omega_X(\cD_X)\otimes_{\cD_X}\cM\hidewidth\cr
\hidewidth\cN\otimes_{\cD_X}\omega_X^{-1}(\cD_X)=\omega^{-1}_X(\cN)
&\longmapsot  &\cN~.\hfill\cr
}
$$

\medskip
{\bf 3.2 De Rham functor for right and left $\cD$-modules}
\smallskip
Let $(\cE,\nabla)$ be a connection on $X$.
By definition its De Rham complex is
$\cOmega^\point_X(\cE)=
\cOmega^\point_X\otimes_{\cO_X}\cE$ where
the differentials are induced by the connection $\nabla$ 
as usual: 
$\nabla(\omega\otimes e)=d(\omega)\otimes e+(-)^{\deg\omega}\omega\wedge\nabla(e)$.
\endgraf
The De Rham functor for left $\cD$-modules
is defined to be compatible with the notion of De Rham complex for
connections, up to a shift.
Let us consider $\cD_X$ as a left $\cD_X$-module,
then its De Rham complex as a connection is
$\cOmega^\point_X(\cD_X)=\cOmega_X^\point\otimes_{\cO_X}\cD_X$
(usual differentials).
It is a resolution of $\omega_X[-\dim X]$ in $\rMod\cD_X$.
For this reason, it is usual  to define
$\DR_X(\cD_X)=\cOmega^\point_X(\cD_X)[\dim X]$,
so that $\DR_X(\cD_X)$ is a resolution of $\omega_X$ in $\rMod\cD_X$.
\endgraf
Now,  if $\cM$ is a left $\cD_X$-module we define
$\DR_X(\cM)=\DR_X(\cD_X)\otimes_{\cD_X}\cM$ which is a complex of
$K_X$-vector spaces.
This functor extends to complexes and
gives in the derived categories the functor
$\DR_X: {\bf D}(\cD_X\lMod) \to {\bf D}(K_X)$
(where ${\bf D}(K_X)$ is the derived category of sheaves in $K_X$-vector
spaces).
Let us observe  that for any $\cM\in {\bf D}(\cD_X\lMod)$ we have
$\DR_X(\cM)=\omega_X\otimes_{\cO_X}^{\bL}\cM$.
\endgraf
The De Rham functor for right $\cD$-modules is defined to be compatible
with the left/right equivalence.
Let us consider $\cD_X$ as a right $\cD_X$-module, then  its De Rham
complex
is $\cTheta^\point_X(\cD_X)=\cD_X\otimes_{\cO_X}\cTheta_X^\point$
(where $\cTheta_X^{\point}=\bigwedge^{-\point}\cDer_X$
and the differentials are locally defined by a Koszul complex).
It is a resolution of $\cO_X$ in $\cD_X\lMod$.
Now, if $\cN$ is a right $\cD_X$-module, then
$\DR_X(\cN)=\cN\otimes_{\cD_X}\cTheta_X^{\point}(\cD_X)$ as a functor
$\rMod\cD_X\to{\bf C}(K_X)$.
The definition naturally extends to the category of complexes of
$\rMod\cD_X$, and to the derived category as before.
\endgraf
The compatibility between De Rham functors is expressed by
the relations
$\DR_X(\cM)=\DR_X(\omega_X(\cM))$
and
$\DR_X(\cN)=\DR_X(\omega_X^{-1}(\cN))$.
\par
{\bf 3.2.1 Relative De Rham functor. }
Let $f:X\to Y$ be a smooth morphism between smooth  $K$-varieties.
The morphism $ f^\ast(\cOmega^1_Y)\rightarrow \cOmega^1_X$
induces a canonical  short exact sequence
$$
0\to f^\ast(\cOmega^1_Y)\to \cOmega^1_X\to \cOmega^1_{X/Y} \to 0
\leqno{\bf (3.2.2)}
$$
where $\cOmega^1_{X/Y}$ is the sheaf of relative differential forms of
degree one.
  Moreover
any $\cO_X$-module in (3.2.2) is locally free of finite type so
(3.2.2) locally splits.
Let $\Theta_{X/Y}=\cHom_{\cO_X}(\cOmega^1_{X/Y},\cO_X)$ be  the
$\cO_X$-dual of
$\cOmega^1_{X/Y}$
and let $\cD_{X/Y}$ denote the $\cO_X$-algebra generated by $\cO_X$ and
$\Theta_{X/Y}$.
As in 1.1 for any $\cO_X$-module $\cE$ the following
supplementary structures on $\cE$ are equivalent:
\item{(1)}
an integrable relative connection, that is a morphism of
$f^{-1}(\cO_Y)$-modules
$\nabla_{X/Y}:{\cE}\to\cOmega^1_{X/Y}\otimes_{{\cO}_X}{\cE}$
Leibniz with respect to sections of ${\cO}_X$ and such that
$\nabla_{X/Y}^2=0$
for the natural extension of $\nabla_{X/Y}$;
\item{(2)}
a structure of left  $\cD_{X/Y}$-module on $\cE$.
\endgraf\noindent
For any $\cO_X$-module $\cE$ endowed with a relative integrable connection
$\nabla_{X/Y}$
let us define its relative
De Rham complex as
$\cOmega^\point_{X/Y}(\cE)=\cOmega^\point_{X/Y}\otimes_{\cO_X}\cE$
where
the differentials are induced by $\nabla_{X/Y}$.
\endgraf
The relative De Rham functor for left $\cD_{X/Y}$-modules
is defined to be compatible with the De Rham functor for connection up to
a shift
(as in the case of $\cD_X$-modules), and it induces a functor of derived
categories:
$$\matrix{
\DR_{X/Y}: & {\bf D}(\cD_{X/Y}) & \longrightarrow& {\bf D}(f^{-1}(\cO_Y))
\hfill\cr
&\hfill \cE & \longmapsto &
\DR_{X/Y}(\cE)=\cOmega^\point_{X/Y}(\cE)[d_{X/Y}]\hfill \cr
}$$
where $d_{X/Y}$ is the relative dimension $d_X-d_Y$.

\medskip
{\bf 3.3 Direct images for connections (the Gauss-Manin connections)}
\medskip
The Leray filtration $\Ler$ on $\cOmega_X^\point$ is defined by
$
\Ler^p\cOmega_X^\point
=\Im(f^\ast\cOmega_Y^p\otimes_{\cO_X}\cOmega_X^{\point-p}\to\cOmega_X^\point)
$
and, since $f$ is a smooth morphism, the associated graded $\cO_X$-module
has
$
\Gr^p_\Ler\cOmega_X^\point\cong
f^\ast\cOmega_Y^p\otimes_{\cO_X}\cOmega_{X/Y}^{\point-p}
$.
\endgraf
If $\cE$ is a connection,  we define the Leray filtration on its
De Rham complex $\cOmega^\point_X(\cE)$  by the tensor product:
$\Ler^p\cOmega^\point_X(\cE)=\Ler^p\cOmega_X^\point\otimes_{\cO_X}\cE$.
Therefore the graded pieces are
$\Gr_\Ler^p\cOmega^\point_X(\cE)=
f^\ast\cOmega_Y^p\otimes_{\cO_X}\cOmega^\point_{X/Y}(\cE)[-p]$.
\endgraf
The Leray filtration induces the
Leray spectral sequence for the direct image functor by $f$:
$$
E_1^{p,q}= \cOmega^p_Y(R^{p+q}f_\ast\cOmega^\point_{X/Y}(\cE)[-p])
\Longrightarrow
\bR^{n}f_\ast\cOmega^\point_{X}(\cE)
\leqno{\bf (3.3.1)}
$$
in the category of $\cO_Y$-modules with differential operators
(the complexes appearing in the spectral sequences
are differential complexes on $Y$).
The differentials of $E_1^{p,q}$ define the Gauss-Manin connections on the
$\cO_Y$-modules $R^qf_\ast\cOmega^\point_{X/Y}(\cE)$, since
$$
d_1^{p,q}:
E_1^{p,q}= \cOmega^p_Y(R^{p+q}f_\ast\cOmega^\point_{X/Y}(\cE)[-p])
\longrightarrow
E_1^{p+1,q}=
\cOmega^{p+1}_Y(R^{p+q+1}f_\ast\cOmega^\point_{X/Y}(\cE)[-p-1])
\leqno{\bf (3.3.2)}
$$
(it is explicitly given by the connection homomorphism
for the direct image functor
of the short exact sequence of complexes
$$
0
\longrightarrow
\Gr_\Ler^{p+1}\cOmega^\point_X(\cE)
\longrightarrow
\Ler^p\cOmega^\point_X(\cE)/\Ler^{p+2}\cOmega^\point_X(\cE)
\longrightarrow
\Gr_\Ler^p\cOmega^\point_X(\cE)
\longrightarrow
0
\leqno{\bf (3.3.3)}
$$
which gives a piece of $E_1$).
\medskip
{\bf 3.4 Direct images for $\cD$-modules}
\medskip
The direct image for $\cD$-modules is defined using the following transfer
modules:
\itemitem{\bf (3.4.1)}
$\cD_{X{\to}Y}=\cO_X\otimes_{f^{-1}\cO_Y}f^{-1}\cD_Y$
which is in $\cD_X\bMod f^{-1}\cD_Y$
(the left $\cD_X$-module structure is induced by (the tensor with) that of
$\cO_X$,
the right $f^{-1}\cD_Y$-module structure is induced by that of
$f^{-1}\cD_Y$,
and the compatibility is obvious); 
\itemitem{\bf (3.4.2)}
$\cD_{Y{\ot}X}=
\omega_X(\cD_X)
\otimes_{\cD_X}
\cD_{X{\to}Y}
\otimes_{f^{-1}\cD_Y}
f^{-1}\omega_Y^{-1}(\cD_Y)$
which is in $f^{-1}\cD_Y\bMod\cD_X$
since it is obtained from $\cD_{X{\to}Y}$ by a double left/right exchange.
\endgraf\noindent
For a right $\cD_X$-module $\cN$, the direct image by $f$ is defined by
$f_+\cN=\bR f_\ast(\cN\otimes_{\cD_X}^\bL\cD_{X{\to}Y})$ as a right
$\cD_Y$-module.
For a left $\cD_X$-module $\cM$, the direct image by $f$ is defined by
$f_+\cM=\bR f_\ast(\cD_{Y{\ot}X}\otimes_{\cD_X}^\bL\cM)$ as a left
$\cD_Y$-module.
The compatibility of the two definitions is expressed by the following
relations:
$$
\omega_Y^{-1}(f_+(\cN))
=f_+(\omega_X^{-1}(\cN))
\quad\hbox{and}\quad
\omega_Y(f_+(\cM))
=f_+(\omega_X(\cM))\ .
$$

\par
\smallskip\eject
{\bf 3.5 Differential complexes and M.Saito equivalence}
\smallskip
Following M. Saito
let $\cO_X\lDiff$ be the category of $\cO_X$-modules with differential
operators as morphisms
and let ${\bf C}(\cO_X\lDiff)$ be its category of complexes (see
[DMSS]).
Objects in ${\bf C}(\cO_X\lDiff)$ are called differential complexes on
$X$.
Any object in  ${\bf C}(\cO_X\lDiff)$ could be regarded as a complex of
$K_X$-vector spaces
so that one has a functor $F:{\bf C}(\cO_X\lDiff)\to{\bf C}(K_X)$.
\endgraf
  In [Sa,1.3.2]    M.Saito defined the linearization functor
$\tDR_X^{-1}:{\bf C}(\cO_X\lDiff)\to{\bf C}(\rMod\cD_X)$ acting on the
differential complex $\cL$
by $\cL\otimes_{\cO_X}\cD_X$ (the differentials being extended to
$\cD_X$-linear maps).
By localization with respect to the multiplicative system of
quasi-isomorphisms on the
right hand side
and with respect to their pull-back 
on the left hand side
(that is, the multiplicative system of $\cD_X$-quasi-isomorphisms:
morphisms of differential complexes whose linearization is a
quasi-isomorphism
in the category of right $\cD_X$-modules)
one obtains the functor
$\tDR_X^{-1}:{\bf D}(\cO_X\lDiff)\to{\bf D}(\rMod\cD_X)$.
\endgraf
It is clear that for any right (resp. left) $\cD_X$-module $\cN$ (resp.
$\cM$),
its De Rham complex is an object in ${\bf C}(\cO_X\lDiff)$.
In particular the functor $\DR_X$ factors through ${\bf C}(\cO_X\lDiff)$
and we denote this factorization by
 $\tDR_X:{\bf C}(\rMod\cD_X)\to{\bf C}(\cO_X\lDiff)$.
In particular we have the following commutative diagram of functors
$$
\xymatrix{
{\bf C}(\rMod\cD_X)\ar[r]^{\tDR_X} \ar[rd]_{\DR_X}&{\bf C}(\cO_X\lDiff)\ar[d]^{F}  \\
&{\bf C}(K_X) \\
}
$$
M.Saito also proved that  $\tDR_X$ localizes with respect to
quasi-isomorphisms
on the left hand side and with respect to $\cD_X$-quasi-isomorphisms on
the right  hand side,
 so it induces
a functor $\tDR_X:{\bf D}(\rMod\cD_X)\to{\bf D}(\cO_X\lDiff)$.

\smallskip
\proclaim{\bf 3.5.1 Proposition (M.Saito's equivalence)}. {
The functor
$\tDR_X^{-1}$ is an equivalence of categories whose quasi-inverse is
$\tDR_X$.
In particular we have canonical quasi-isomorphisms
$$
 \tDR_X^{-1}\tDR_X(\cN)\cong\cN
\qquad\hbox{and}\qquad
\tDR_X^{-1}\tDR_X(\cM)\cong\omega_X(\cM)\ ,
$$
both  in ${\bf C}(\rMod\cD_X)$
(see [Sa,1.8]).
}\par
\smallskip

{\bf 3.5.2 }
In the case of right $\cD_X$-modules, there is also a compatibility with
the
direct image of differential complexes
(which is induced by the usual direct image for abelian sheaves),
via the linearization functor:
$\tDR_Y^{-1}\bR f_\ast(\cL)=f_+\tDR_X^{-1}(\cL)$
(see [DMSS,1.3.2]).

\smallskip
\proclaim{\bf 3.6. Theorem (comparison for direct images, following
[DMSS])}. {
Let $f:X\to Y$ be a smooth morphism of smooth $K$-varieties.
For any left $\cD_X$-module $\cM$ (identified with a connection on $X$)
and for any $q$ we have natural isomorphisms
$
R^{q}f_\ast\DR_{X/Y}(\cM)\cong\cH^q(f_+\cM)
$
in the category of left $\cD_Y$-modules,
where the left hand side has the structure of the Gauss-Manin connection.
}\par
{\bf Proof.}
Let consider the Leray spectral sequence $E$ of $\cM$ with respect to $f$.
Since $\cD_Y$ is a flat $\cO_Y$-module,
we may apply the linearization functor $\tDR_Y^{-1}$ to obtain the spectral
sequence
$$
\tDR_Y^{-1}E_1^{p,q}=
\tDR_Y^{-1}\cOmega^p_Y(R^{p+q}f_\ast\cOmega^\point_{X/Y}(\cM)[-p])
\Longrightarrow
\tDR_Y^{-1}R^{n}f_\ast\cOmega^\point_{X}(\cM)
$$
in the category of right $\cD_Y$-modules.
Now (by 3.5.1) the complex $\tDR_Y^{-1}E_1^{\point,q}$ is quasi-isomorphic to
$$
\tDR_Y^{-1}\tDR_Y[-\dim Y](R^qf_\ast\DR_{X/Y}[-d_{X/Y}](\cM))
\cong
\omega_Y(R^qf_\ast\DR_{X/Y}(\cM))[-\dim Y-d_{X/Y}]
$$
(where $d_{X/Y}=\dim X-\dim Y$ is the relative dimension)
so that the spectral sequence degenerates at $E_2$;
while the limit is quasi-isomorphic (by 3.5.2, 3.5.1 and 3.4) to
$$
\eqalign{
\tDR_Y^{-1}\bR f_\ast\tDR_{X}[-\dim X](\cM)
&\cong
f_+(\tDR_X^{-1}\tDR_{X}[-\dim X](\cM))
\cr &\cong
f_+( \omega_X(\cM)[-\dim X])
\cr &\cong
\omega_Y(f_+(\cM)[-\dim X]) \ .
}$$
So we have the isomorphisms
$
R^qf_\ast\DR_{X/Y}(\cM)
\cong
\cH^{q+\dim X}(f_+(\cM)[-\dim X])
$
in the category of complexes of left $\cD_Y$-modules,
from which the proposition follows.
\hfill$\square$

\bigskip
{\bf\S 4. Alternative proof of the comparison between direct images}
\bigskip

We present here an alternative proof of the comparison theorem for direct
images, which is
in some sense more elementary, since Saito's equivalence is not
used. The strategy we discuss here also has the advantage of clarifying
the structure of the Gauss-Manin connection,
taking account of one of its avatars before the application of the derived
direct image functor.
In fact we compare two structures of left $f^{-1}\cD_Y$-modules,
one defining the Gauss-Manin connection,
the other coming from the structure of the transfer module $\cD_{Y\ot
X}$.
The main technical tool is  the commutativity of a diagram
in the derived category of right $\cD_X$-modules 
(see Proposition 4.3.2), for which the homotopy lemma 4.3.1 is
used.
\endgraf

\medskip
{\bf 4.1 The distinguished triangle for the Gauss-Manin connection}
\medskip
Let us consider the exact sequence (3.3.3) for $p=0$
and $\cE=\cD_X$ (in the category of complexes of right $\cD_X$-modules):
$$
0
\longrightarrow
f^\ast(\cOmega^1_Y)\otimes_{\cO_X}\cOmega^\point_{X/Y}(\cD_X)[-1]
\mathop{\longrightarrow}\limits^i
\Ler^0\cOmega^\point_X(\cD_X)/\Ler^{2}\cOmega^\point_X(\cD_X)
\mathop{\longrightarrow}\limits^\pi
\cOmega^\point_{X/Y}(\cD_X)
\longrightarrow
0
\leqno{\bf (4.1.0)}
$$
(recall that
$\Gr_\Ler^1(\cOmega^\point_{X}(\cD_X))\isom
f^\ast(\cOmega^1_Y)\otimes_{\cO_X}\cOmega^\point_{X/Y}(\cD_X)[-1]$
and
$\Gr_\Ler^0(\cOmega^\point_{X}(\cD_X))\isom\cOmega^\point_{X/Y}(\cD_X)$)
which gives the distinguished triangle:
$$
\xymatrix{
f^\ast(\cOmega^1_Y)\otimes_{\cO_X}\cOmega^\point_{X/Y}(\cD_X)[-1]
\ar[r]^i & \Ler^0\cOmega^\point_X(\cD_X)/\Ler^{2}\cOmega^\point_X(\cD_X)
\ar[d]^\pi \\
& \cOmega^\point_{X/Y}(\cD_X) \ar[lu]_{[+1]}^{\delta(\cD_X)} \\}
\leqno{\bf (4.1.1)}
$$
in ${\bf D}^b(\rMod\cD_X)$
(derived category of right $\cD_X$-modules).
The connecting morphism $\delta(\cD_X)$ is represented in the derived
category
by the following diagram (in which the mapping cone of $i$ appears):
$$
\xymatrix{
(f^\ast(\cOmega^1_Y)\otimes_{\cO_X}\cOmega^\point_{X/Y}(\cD_X))
\oplus
\Ler^0\cOmega^\point_X(\cD_X)/\Ler^{2}\cOmega^\point_X(\cD_X)
\ar[d]_{\rm qis}^{(0,\pi)} \ar[rd]^{-p_1=(-1,0)}
\\
\cOmega^\point_{X/Y}(\cD_X)&
f^\ast(\cOmega^1_Y)\otimes_{\cO_X}\cOmega^\point_{X/Y}(\cD_X)\\}
\leqno{\bf (4.1.2)}
$$
The fundamental fact here is that
{\sl the Gauss-Manin connection comes by
applying the derived functor $\bR f_\ast$ to this connecting morphism. }
\endgraf
Notice moreover that for any left $\cD_X$-module $\cE$ we can apply the
derived
functor
${\rm -}\otimes_{\cD_X}^{\bL}\cE$ to the distinguished triangle (4.1.1)
so we obtain the distinguished triangle:
$$
\xymatrix{
f^\ast(\cOmega^1_Y)\otimes_{\cO_X}\cOmega^\point_{X/Y}(\cE)[-1]
\ar[r] & \Ler^0\cOmega^\point_X(\cE)/\Ler^{2}\cOmega^\point_X(\cE) \ar[d]
\\
& \cOmega^\point_{X/Y}(\cE) \ar[lu]_{[+1]}^{\delta(\cE)} \\}
$$
(really, we do not need to take the derived functor ${\rm
-}\otimes_{\cD_X}^{\bL}\cE$
but simply ${\rm -}\otimes_{\cD_X}\cE$ because any complex
in (4.1.1) is acyclic for ${\rm -}\otimes_{\cD_X}\cE$).
This is the distinguished triangle induced by the short exact sequence
(3.3.3)
for $p=0$.

\medskip
{\bf 4.2 The connection associated to $\cD_{Y{\ot}X}$ }
\medskip
In view of the comparison, we have to give a concrete expression in terms
of a connection (in the derived category) for the $\cD_Y$-module structure
of the derived direct image.

\proclaim{\bf 4.2.1 Lemma [Lau, 5.2.3.4]}. {
Let $f:X\to Y$ be a smooth morphism of smooth $K$-varieties.
There exists a canonical morphism of complexes of right $\cD_X$-modules
$\lambda:\DR_{X/Y}(\cD_X)\rightarrow \cD_{Y{\ot}X}$
which is a quasi-isomorphism.
In particular, $\DR_{X/Y}(\cD_X)$ is a left resolution of
$\cD_{Y{\ot}X}$
in the category $f^{-1}\cO_Y\bMod\cD_X$,
with locally free right $\cD_X$-modules.
As a consequence, the complex $\DR_{X/Y}(\cD_X)$ admits a structure of
left $f^{-1}\cD_Y$-module,
induced by transfer of structure via the quasi-isomorphism,
in the derived category of $f^{-1}\cO_Y$-Mod.
}\par
{\bf Proof.}
The canonical morphism $\lambda$ is defined by the composition 
$$
\xymatrix{
\omega_{X/Y}\otimes_{\cO_X}\cD_X 
 \ar[r]^(.45){i} 
&
\cD_{Y{\ot}X}\otimes_{\cO_X}\cD_X
 \ar[r]^(.6){\cdot} 
&
\cD_{Y{\ot}X}\\
}
$$
where the first map comes from the canonical inclusion of $\omega_{X/Y}$ 
into $\cD_{Y{\ot}X}=\omega_{X/Y}\otimes_{\cO_X}f^\ast\cD_Y$, 
and the second one uses 
the canonical structure of right $\cD_X$-module of 
$\cD_{Y{\ot}X}$.
A local computation using the canonical filtrations by the order of
differential operators
shows that the graded pieces are Koszul complexes, so that the assertion
follows
(see [Lau, 5.2.3.4] for details).
\hfill$\square$ 
\smallskip
Notice moreover that the morphism $\lambda$ may be described 
in terms of the canonical map $\cD_f:\cD_X\to f^\ast\cD_Y$ 
composing $\id\otimes\cD_f$ with the transposition of 
differential operators, since the following diagram 
\def\struttura{{{\scriptstyle\blacktriangledown}\atop ~}}
$$
\xymatrix{ 
\omega_{X/Y}\otimes_{}\cD_X \isom
\cHom_{}(f^{-1}\omega_Y,\omega_{X}\mathop\otimes\limits^{\struttura}_{}\cD_X) 
 \ar[r]^(.45){i} \ar[d]_{\iota}
&
\cHom_{}(f^{-1}\omega_Y,(\omega_{X}\otimes f^{\ast}\cD_Y^{\struttura})\otimes_{}\cD_X) 
 \ar[d]^{\cdot} 
\isom\cD_{Y{\ot}X}\otimes_{}\cD_X
\\ 
\cHom_{}(f^{-1}\omega_Y,\omega_{X}\otimes_{}\cD_X^{\struttura}) 
 \ar[r]_{\cD_f} 
&\cHom_{}(f^{-1}\omega_Y,\omega_{X}\otimes f^{\ast}\cD_Y^{\struttura}) 
\isom\cD_{Y{\ot}X}
\\
}
$$
commutes (here, $\cHom$ means $\cHom_{f^{-1}\cO_Y}$, 
the symbol ${}^\struttura$ indicates the $f^{-1}\cO_Y$-module 
structure used in the $\cHom$, 
$\otimes$ means $\otimes_{\cO_X}$,
and 
$\iota$ is the canonical involution of 
$\omega_{X}\otimes_{\cO_X}\cD_X$ exchanging the two structures of 
right $\cD_X$-module). 

\smallskip
{\bf 4.2.2 }
Notice that the  $f^{-1}\cD_Y$-module structure
of $\DR_{X/Y}(\cD_X)$ is described in terms of its connection by the
following diagram
$$
\xymatrix{
\DR_{X/Y}(\cD_X) \ar[d]_{\rm qis}^{\lambda}
&
f^{-1}\cOmega^1_Y\otimes_{f^{-1}\cO_Y} \DR_{X/Y}(\cD_X)
\ar[d]_{\rm qis}^{1\otimes\lambda}\\
\cD_{Y{\ot}X} \ar[r]^(.3){\nabla_{Y{\ot}X}}
&
f^{-1}\cOmega^1_Y\otimes_{f^{-1}\cO_Y} \cD_{Y{\ot}X}\\
}
$$
where $\nabla_{Y{\ot}X}$ is the connection
(Leibniz with respect to section of $f^{-1}(\cO_Y)$)
induced by the $f^{-1}(\cD_Y)$-module structure of $\cD_{Y\ot X}$.
\endgraf 

\smallskip
\proclaim{\bf 4.2.3 Corollary}. {
Let $f:X\to Y$ be a smooth morphism of smooth $K$-varieties.
For any left $\cD_X$-module $\cM$ there is a canonical quasi-isomorphism
$\lambda(\cM):\DR_{X/Y}(\cM)\to\cD_{Y{\ot}X}\otimes_{\cD_X}^\bL\cM$,
so that $\DR_{X/Y}(\cM)$ admits a structure of left $f^{-1}(\cD_Y)$-module
in the derived category of $f^{-1}(\cO_Y)$-modules.
\endgraf
Applying the derived direct image functor to the above morphism,
one obtains  a canonical
quasi-isomorphism
$\bR f_\ast\DR_{X/Y}(\cM)\to f_{+}(\cM)$ in the category of
complexes of $\cO_Y$-modules,
so that the complex $\bR f_\ast\DR_{X/Y}(\cM)$ admits
a structure of left $\cD_Y$-module in the derived category of
$\cO_Y$-modules.
}\par

\medskip
{\bf 4.3 }
The following two results are the kernel of our comparison argument;
essentially we have to compare the diagram {(4.1.2)}
with the diagram in 4.2.2.
To do that, we re-write the projection $p_1$ of {(4.1.2)} up to homotopy.

\proclaim{\bf 4.3.1 (homotopy) Lemma}. {
Let
$i:f^\ast(\cOmega^1_Y)\otimes_{\cO_X}\cOmega^\point_{X/Y}(\cD_X)[-1]
\rightarrow
\Ler^0\cOmega^\point_X(\cD_X)/\Ler^{2}\cOmega^\point_X(\cD_X)$
be the canonical inclusion.
\item{(i)}
the identity morphism of the mapping cone of $i$
is homotopic to the morphism $\Psi^\point$ defined by
$\Psi^{d_{X/Y}}=\left({0\atop 0}{-\phi^{-1}d\atop 1}\right)$
and $\Psi^{q}=\left({1\atop 0}{0\atop 1}\right)$
for $q\neq d_{X/Y}$;
\item{(ii)}
the connecting morphism
$$
-p_1=(-1,0):
(f^\ast(\cOmega^1_Y)\otimes_{\cO_X}\cOmega^\point_{X/Y}(\cD_X))
\oplus
\Ler^0\cOmega^\point_X(\cD_X)/\Ler^{2}\cOmega^\point_X(\cD_X)
\rightarrow
(f^\ast(\cOmega^1_Y)\otimes_{\cO_X}\cOmega^\point_{X/Y}(\cD_X))
$$
of the distinguished triangle generated  by $i$
is homotopic to the morphism $\psi^\point=-p_1\circ\Psi^\point$,
that is
$\psi^{d_{X/Y}}=(0,\phi^{-1}d)$
and $\psi^{q}=(-1, 0)$
for $q\neq d_{X/Y}$;
\endgraf\noindent
where $\phi$ is the canonical isomorphism
$
f^\ast\cOmega^1_Y\otimes_{\cO_X}\omega_{X/Y}(\cD_X)
\longrightarrow
{\Ler^0\over \Ler^2}(\cOmega^{d_{X/Y}+1}_X(\cD_X))
$
induced by $i$.
}\par
{\bf Proof.}
Let us consider the exact sequence of complexes (4.1.0)
in degrees $d_{X/Y}-1, d_{X/Y}, d_{X/Y}+1$:
$$
\xymatrix{
0\ar[r] &
f^\ast\cOmega^1_Y\otimes_{\cO_X}\cOmega^{d_{X/Y}-2}_{X/Y}(\cD_X)
\ar[r]^i\ar[d]^d &
{\Ler^0\over \Ler^2}(\cOmega^{d_{X/Y}-1}_{X}(\cD_X)) \ar[r]^\pi \ar[d]^d &
 \cOmega_{X/Y}^{d_{X/Y}-1}(\cD_X) \ar[r]\ar[d]^d &
 0 \\
0\ar[r] &
f^\ast\cOmega^1_Y\otimes_{\cO_X}\cOmega^{d_{X/Y}-1}_{X/Y}(\cD_X)
\ar[r]^i\ar[d]^d &
{\Ler^0\over \Ler^2}(\cOmega^{d_{X/Y}}_{X}(\cD_X)) \ar[r]^\pi \ar[d]^d &
 \omega_{X/Y}(\cD_X) \ar[r]\ar[d] &
 0 \\
0\ar[r] &
f^\ast\cOmega^1_Y\otimes_{\cO_X} \omega_{X/Y}(\cD_X)
\ar[r]_(0.5){\cong}^(.5){\phi}  &
{\Ler^0\over \Ler^2}(\cOmega^{d_{X/Y}+1}_{X}(\cD_X))  \ar[r]  &
0 \ar[r]&
0 \\
}
$$
from which one deduces that $\phi$ is an isomorphism.
\item{\sl (i)}
Using the homotopy map of the mapping cone of $i$ which is zero for
degrees
different from $d_{X/Y}+1$, and ${-\phi^{-1}\choose 0}$ in degree
$d_{X/Y}+1$,
we have that the identity map of the mapping cone is homotopic to
the morphism having the following expression in degree $d_{X/Y}$:
$
\left({1\atop 0}{0\atop 1}\right)+{-\phi^{-1}\choose 0}(\phi, d)=
\left({1\atop 0}{0\atop 1}\right)+\left({-1\atop 0}{-\phi^{-1}d\atop
0}\right)=
\left({0\atop 0}{-\phi^{-1}d\atop 1}\right)
$
and unchanged otherwise,
as stated.
\item{\sl (ii)}
This follows from the previous point, since $\id\sim\Psi^\point$
implies $p_1\sim p_1\circ\Psi^\point =\psi^\point$.
Explicitly,  we may use the homotopy map of the
connecting morphism which is zero for degrees
different from $d_{X/Y}+1$, and $\phi^{-1}$ in degree $d_{X/Y}+1$.
Then we have that $-p_1$ is homotopic to
the morphism having the following expression in degree $d_{X/Y}$:
$
(-1,0)+(\phi^{-1})(\phi,d)=(0,\phi^{-1}d)
$
and unchanged otherwise,
as stated.
\hfill$\square$

\proclaim{\bf 4.3.2 Proposition}. {
The following diagram
$$
\xymatrix{
(f^\ast(\cOmega^1_Y)\otimes_{\cO_X}\DR_{X/Y}(\cD_X))
\oplus
\Ler^0\DR_X(\cD_X)/\Ler^{2}\DR_X(\cD_X)[-d_Y]
\ar[d]_{\rm qis}^{(0,\pi)} \ar[rd]^{-p_1=(-1,0)}
\\
\DR_{X/Y}(\cD_X) \ar[d]_{\rm qis}^{\lambda}
&
f^\ast(\cOmega^1_Y)\otimes_{\cO_X} \DR_{X/Y}(\cD_X)
\ar[d]_{\rm qis}^{1\otimes\lambda}\\
\cD_{Y{\ot}X} \ar[r]^(.3){\nabla_{Y{\ot}X}}
&
f^\ast(\cOmega^1_Y)\otimes_{\cO_X}\cD_{Y{\ot}X}\\
}
$$
commutes in the derived category of right $\cD_X$-modules.
}\par

Since this diagram is just the superposition of (4.1.2) and 4.2.2,
its commutativity identifies the two structures of left $f^{-1}\cD_Y$
module of $\DR_{X/Y}(\cD_X)$.
Moreover, the result extends immediately to any left $\cD_X$-module
$\cM$, by the tensor product over $\cD_X$.

{\bf Proof.}
By the homotopy lemma, we may use $\psi$ instead of $-p_1$,
and since any object of the last row is a complex concentrated in degree
zero,
we need only  prove the commutativity of the following diagram:
$$
\xymatrix{
{\Ler^0\over \Ler^2}(\cOmega^{d_{X/Y}}_{X}(\cD_X))
\ar[d]^{\pi} \ar[rd]^{\phi^{-1}d}
\\
\omega_{X/Y}(\cD_X) \ar[d]^{\lambda}
&
f^\ast(\cOmega^1_Y)\otimes_{\cO_X} \omega_{X/Y}(\cD_X)
\ar[d]^{1\otimes\lambda}\\
\cD_{Y{\ot}X} \ar[r]^(.3){\nabla_{Y{\ot}X}}
&
f^\ast(\cOmega^1_Y)\otimes_{\cO_X} \cD_{Y{\ot}X}\\
}
$$ 
(the first object of the mapping cone does not appear, 
because it is sent to zero using both morphisms). 
\endgraf 
Notice, first of all, that there exists a unique morphism
$\alpha:\omega_{X/Y}(\cD_X)\to f^{\ast}\cOmega^1_Y\otimes_{\cO_X}
\cD_{Y{\ot}X}$
making the upper part (``parallelogram'') of the diagram commutative,
since
(with reference to the diagram in 4.3.1) we have that
$(1\otimes\lambda)\phi^{-1}di=(1\otimes\lambda)d=0$,
so that $(1\otimes\lambda)\phi^{-1}d$ factors  as $\alpha\pi$.
As a consequence,
we have to prove that $\nabla_{Y{\ot}X}\lambda=\alpha$.
One has two different possibilities  for the proof. The first one is a
local computation. 
Using local coordinates $x_i$ ($i=1,\dots,d_X$) 
on $X$ such that $dx_1,\dots,dx_{d_{X/Y}}$ 
are generators of the relative differentials, 
and $\omega=dx_1\wedge\cdots\wedge dx_{d_{X/Y}}$, 
for the morphism $\alpha$ we have the following expression 
$$
\eqalign{\alpha(\omega\otimes\partial)
&= (1\otimes\lambda)\phi^{-1}d([\omega\otimes\partial]) \cr
&=(1\otimes\lambda)\phi^{-1}\left[
\sum\nolimits_idy_i\wedge\omega\otimes\eta_i\partial
\right]\cr
&=(1\otimes\lambda)\left(
\sum\nolimits_idy_i\otimes\omega\otimes\eta_i\partial
\right) \cr
&=\sum\nolimits_idy_i\otimes\omega\otimes(\eta_i\partial)^\ast
}$$
where local coordinates $y_i$ on $Y$ are used 
($dy_i$ and $\eta_i$ are dual bases for $\cOmega^1_Y$ and$\cTheta_Y$); 
and $\partial$ is a local section of $\cD_X$. 
Therefore, 
the action of the derivative $\eta_i$ (dual bases of $dy_i$) 
is given by 
$
\omega\otimes(\eta_i\partial)^\ast$. 
On the other hand, the structure of $f^{-1}\cD_Y$-module of 
$\cD_{Y{\ot}X}$ 
is given by the twist of the right structure (by multiplication) of 
$f^\ast\cD_Y$ 
with the right structure of $f^{-1}\omega_Y$ 
(by the action of $-{\rm Lie}_\eta$, which is trivial on $\omega$), 
and has therefore the local expression 
$\eta_i(\omega\otimes\partial)=
-\omega\otimes\partial\eta_i$; 
composing with $\lambda$, the action of $\eta_i$ sends 
$\omega\otimes\partial$ to 
$\omega\otimes\partial^\ast\eta_i^\ast$
and coincides with that given above.
\endgraf
The second possibility is more abstract and relies on the compatibility of
De Rham functors.
Since the morphism $\alpha$ is induced by the differential of the absolute
De Rham complex
of $\cD_X$, and it factors uniquely through $\lambda$ 
(see again the diagram in 4.3.1, since $(\pi)d\alpha=0$ 
and $\lambda$ gives a cokernel for $d$ in the last column),
it is sufficient to prove that the morphism $\nabla_{Y{\ot}X}$ is
compatible with that complex.
We observe that $\omega_X(f^\ast\cD_Y)$
has two compatible structures,
as a right  $\cD_X$-module and as a right $f^{-1}(\cD_Y)$-module.
We define
$\DR_Y(\omega_X(f^\ast\cD_Y)):=
\omega_X(f^\ast\cD_Y)\otimes_{\cO_X}\bigwedge^{-\point}f^{\ast}\cTheta^1_Y$ 
(the De Rham complex as right $f^{-1}\cD_Y$-module)
which is isomorphic  (by left/right exchange on the $f^{-1}(\cD_Y)$-module
structure) to
$\DR_{Y}(\omega_{X/Y}(f^\ast\cD_Y)):=
f^{-1}\cOmega^{\point}_Y\otimes_{f^{-1}(\cO_Y)}\omega_{X/Y}(f^\ast\cD_Y)[d_Y]$.
Now, we have the following canonical morphisms of De Rham
complexes:
$$
\DR_X(\cD_X)
\isom
\DR_X(\omega_X(\cD_X))
\mathop{\longrightarrow}\limits^\iota_\sim
\DR_X(\omega_X(\cD_X))
\longrightarrow
\DR_{Y}(\omega_X(f^\ast\cD_Y))
\isom
\DR_{Y}(\omega_{X/Y}(f^\ast\cD_Y))
$$
where the isomorphisms $\isom$ are given by left/right exchanges,
the first morphism comes from the involution $\iota$, 
the second comes from  $\cD_f:\cD_X\to f^\ast\cD_Y$.
In  degrees $-d_Y$ and $-d_{Y}+1$ we can read the following
compatibilities:
\catcode`\@=11
\def\mmatrix#1{\null\,\vcenter{\normalbaselines\m@th
                        \ialign{\hfil$##$\hfil&&\ \hfil$##$\hfil\crcr
                        \mathstrut\crcr\noalign{\kern-\baselineskip }
#1\crcr\mathstrut\crcr\noalign{\kern-\baselineskip}}}\,}
\catcode`\@=12
\def\dmatrix{%
\def\normalbaselines{\baselineskip15pt\lineskip3pt\lineskiplimit3pt}%
                  \mmatrix
}
$$ \hskip-20pt
\dmatrix{
\cOmega_X^{d_{X/Y}}(\cD_X)
&\mathop{\rightarrow}\limits_\sim
&\omega_X(\cD_X)\otimes\bigwedge^{d_Y}\cTheta^1_X
&\to
&\omega_X(f^\ast\cD_Y)\otimes\bigwedge^{d_Y}f^{\ast}\cTheta^1_Y
&\isom
&\omega_{X/Y}(f^\ast\cD_Y)
\cr
\big\downarrow
&&\big\downarrow
&&\big\downarrow
&&\big\downarrow
 \hbox to0pt{$\scriptstyle\nabla_{Y\ot X}\hss$}
\cr
\cOmega_X^{d_{X/Y}+1}(\cD_X)
&\mathop{\rightarrow}\limits_\sim
&\omega_X(\cD_X)\otimes\bigwedge^{d_Y-1}\cTheta^1_X
&\to
&\omega_X(f^\ast\cD_Y)\otimes\bigwedge^{d_Y-1}f^{\ast}\cTheta^1_Y
&\isom
&f^{-1}\cOmega_Y^1\otimes_{}\omega_{X/Y}(f^\ast\cD_Y)
\cr
}
$$
from which the result follows.
\hfill$\square$

\proclaim{\bf 4.4 Theorem}. {
For any left $\cD_X$-module $\cM$,
the canonical quasi-isomorphism
$$\bR f_\ast\DR_{X/Y}(\cM)\to%
\bR
f_\ast(\cD_{Y{\ot}X}\otimes^{\bL}_{\cD_X}\cM)=
f_+(\cM)$$
identifies the structures of $\cD_Y$-modules of the two terms
(the Gauss-Manin structure on the left hand side, and the canonical one on
the
right hand side),
so that for any $q$ we have natural isomorphisms
$$
R^{q}f_\ast\DR_{X/Y}(\cM)\cong\cH^q(f_+\cM)
$$
in the category of left $\cD_Y$-modules,
where the left side has the structure of the Gauss-Manin connection.
}\par
{\bf Proof.}
We have (by lemma 4.2.3) a canonical quasi-isomorphism
$\DR_{X/Y}(\cM)\to\cD_{Y{\ot}X}\otimes^{\bL}_{\cD_X}\cM$
which identifies (by Proposition 4.3.2)
the structures of $f^{-1}\cD_Y$-modules of the two terms.
Applying the derived functor $\bR f_\ast$ one obtains
a canonical quasi-isomorphism
$\bR f_\ast\DR_{X/Y}(\cM)\to\bR
f_\ast(\cD_{Y{\ot}X}\otimes^{\bL}_{\cD_X}\cM)=
f_+(\cM)$
which identifies the structures of $\cD_Y$-modules of the two terms.
(Recall that $\bR f_\ast$ is intended as a functor from
${\bf D}^+(f^{-1}\cD_Y\lMod)$ to ${\bf D}^+(\cD_Y\lMod)$.)
\hfill$\square$

\medskip
{\bf 4.5}
Finally, we translate the result in term of De Rham complexes.

\proclaim{\bf 4.5.1 Lemma [Me ch.I, 5.4.3, 4, 4b]}. {
Let $f:X\to Y$ be a smooth morphism of smooth $K$-varieties.
For any $\cD_X$-module $\cM$ there exists a canonical quasi-isomorphism
$\DR_Y(f_+\cM)\rightarrow\bR f_\ast(\DR_X(\cM))$,
where $\bR f_\ast$ is the usual derived direct image for abelian sheaves.
}\par
{\bf Proof.}
The canonical morphism in the derived category of abelian sheaves on $Y$
$$
\omega_Y\otimes^\bL_{\cD_Y}\bR f_\ast(\cD_{Y{\ot}X}\otimes^\bL_{\cD_X}\cM)
\longrightarrow
\bR
f_\ast(f^{-1}\omega_Y\otimes^\bL_{f^{-1}\cD_Y}(\cD_{Y{\ot}X}\otimes^\bL_{\cD_X}\cM))
$$
is an isomorphism: this can be verified using a locally free resolution of
$\omega_Y$ in the category of right-$\cD_Y$-modules.
Now  the left hand side is canonically quasi-isomorphic to
$\omega_Y\otimes^\bL_{\cD_Y}f_+(\cM)\isom\DR_Y(f_+\cM)$
and the right hand side to
$\bR f_\ast(\omega_X\otimes^\bL_{\cD_X}\cM)\isom\bR f_\ast(\DR_X(\cM))$.
\hfill$\square$
\proclaim{\bf 4.5.2 Corollary}. {
Let $f:X\to Y$ be a smooth morphism of smooth $K$-varieties.
For any $\cD_X$-module $\cM$ there exists a canonical quasi-isomorphism
$
\DR_Y(\bR f_\ast\DR_{X/Y}(\cM))\isom
\DR_Y(f_+\cM)\isom\bR f_\ast(\DR_X(\cM))
$,
where $\bR f_\ast$ is the usual derived direct image for abelian sheaves.
}\par
\medskip
{\bf 4.6 A final remark. } 
A naive approach to the comparison problem 
for direct images 
may be the following argument,
which is very elementary.
The existence of the morphisms of the Theorem 4.4
in the category of $\cO_Y$-modules
and the fact that they are isomorphisms follow from
the Corollary 4.2.3.
In order to conclude the proof, it is enough to identify the structure of
$\cD_Y$-module
induced by the structure of the Gauss-Manin connection on the left
with that of the right, and this may be done locally.
The isomorphisms of  Lemma 4.5.1 can be constructed locally in the
following way.
Using a local decomposition
$f^\ast\cOmega^1_Y\oplus\cOmega^1_{X/Y}\isom\cOmega^1_X$,
we have (locally) that
$(f^\ast\cOmega^\point_Y\otimes_{\cO_X}\cOmega^\point_{X/Y})_\tot
\isom\cOmega^\point_X$
(see [B.A ch.III \S 7 n.7] as graded modules, and the usual signs
convention
identifies the differentials).
In this isomorphism the Leray filtration on the right is identified
with the
filtration $F_I$ of the total complex on the left.
Therefore the Leray spectral sequence of $\cM$ for $f$
is locally identified with the spectral sequence $E_I$ of the bicomplex
$(f^\ast\cOmega^\point_Y\otimes_{\cO_X} \cOmega^\point_{X/Y}(\cM))$
with respect to the functor $\bR f_\ast$.
In particular we have at the $E_1$ level that
$\DR_Y(R^{q}f_\ast\DR_{X/Y}(\cM))\isom\DR_Y(\cH^q(f_+\cM))$,
which identifies the Gauss-Manin connection with the connection induced on
$\cH^q(f_+\cM)$
by the $\cD_Y$-module structure of $f_+\cM$,
and concludes the ``proof''.
The unpleasant point in this argument, 
and probably the reason of its absence in the literature,  is that
it makes use of a non canonical local decomposition of the
canonical exact sequence of differentials.

\bigskip
{\bf References}
\medskip

\item{[B.A]} Bourbaki, Alg\`ebre.
\item{[B.AC]} Bourbaki, Alg\`ebre Commutative.
\item{[Be]} Bernstein,
Algebraic theory of $\cD$-modules, preprint.
\item{[BO]} Berthelot, Ogus, Notes on
Crystalline cohomology.
Mathematical Notes Princeton University Press, 1978.
\item{[Bo]} Borel,
Algebraic $D$-modules.
Perspectives in Mathematics, 2. Academic Press, Inc., Boston, MA,  1987.
xii+355 pp.
\item{[D]} Deligne,
\'Equations diff\'erentielles \`a points singuliers r\'eguliers.
Lecture Notes in Mathematics, Vol. 163. Springer-Verlag, Berlin-New York,
1970. iii+133 pp.
\item{[DMSS]} Dimca, Maaref, Sabbah, Saito,
Dwork cohomology and algebraic $\cD$-modules,
Math. Ann.  318  (2000), 107--125.
\item{[EGA]}
Grothendieck, El\'ements de G\'eometrie Alg\`ebrique,
$I$(new) Springer,
$I\!I$, $I\!I\!I$, $I\!V$ Publ.Math.IHES.
\item{[Lau]} Laumon,
Sur la cat\'egorie d\'eriv\'ee des ${\cal D}$-modules filtr\'es.
Algebraic geometry (Tokyo/Kyoto, 1982),   151--237,
Lecture Notes in Math., 1016, Springer, Berlin,  1983.
\item{[Me]} Mebkhout,
Le formalisme des six op\'erations de Grothendieck pour les $\cD\sb
X$-modules coh\'erents.
Travaux en Cours, 35.
Hermann, Paris,  1989. x+254 pp.
\item{[Sa]} Saito M.,
Induced $\cD$-modules and differential complexes.
Bull. Soc. Math. France  117  (1989),  361--387.

\end